\documentclass[12pt]{amsart}
\usepackage{amsmath}
\usepackage{amsfonts,amssymb,amsthm, txfonts, pxfonts}
\def\struckint{\mathop{%
\def\mathpalette##1##2{\mathchoice{##1\displaystyle##2}%
  {##1\textstyle##2}{##1\scriptstyle##2}{##1\scriptscriptstyle##2}}%
\mathpalette
{\vbox\bgroup\baselineskip0pt\lineskiplimit-1000pt\lineskip-1000pt
\halign\bgroup\hfill$}
{##$\hfill\cr{\intop}\cr\diagup\cr\egroup\egroup}%
}\limits}

\newcommand{\Nmin}{\mathrm{Nmin\:}}

\newtheorem{theorem}{Theorem}
\newtheorem{lemma}[theorem]{Lemma}

\newtheorem{corollary}[theorem]{Corollary}

\newtheorem{definition}[theorem]{Definition}

\theoremstyle{remark}

\newtheorem{remark}[theorem]{Remark}

\newcommand{\vol}{{\mathrm{vol}\,}}

\begin{document}

%-------------- Author entries --------------------

\title{Surface area of Ellipsoids}
%Article title
%\shorttitle{DRAFT} % Shortened version for
                                             % headline title

\author{Igor Rivin}

%\address{Mathematics Department, California Institute of Technology,
%Pasadena, CA 91125}

\address{Department of Mathematics, Temple University, Philadelphia}

\email{rivin@math.temple.edu}

\thanks{The author is supported by the NSF DMS. Parts of this paper appeared
  as the preprint \cite{riv1}; the author would like to thank Warren~D.~Smith
  on comments on a previous version of this paper.The author would also like
  to thank Princeton University for its hospitality.}

\date{\today}



\begin{abstract}
We study the surface area of an ellipsoid in $\mathbb{E}^n$ as the function of
the lengths of their major semi-axes.  We write down an explicit formula as an
integral over $\mathbb{S}^{n-1},$ use this formula to derive convexity
properties of the surface area, to greatly sharpen the estimates given
in \cite{riv2} for the surface area of a large-dimensional ellipsoid,
to produce asymptotic formulas for the surface area and the
\emph{isoperimetric ratio} of an ellipsoid in large dimensions, and to give an
expression for the surface in terms of the Lauricella hypergeometric
function.
\end{abstract}

\maketitle

\section*{Introduction}
In \cite{riv2} estimates were given for the \emph{mean curvature integrals} of
an ellipsoid $E$ in $\mathbb{E}^n$ in terms of the lengths of its \textit{major
  semiaxes} -- the $0$-th mean curvature integral is simply the surface area
of the ellipsoid $E.$ The estimate for the surface area was simply the
$n-1$-th symmetric function of the semi-axes, and it was shown in \cite{riv2}
that this estimate differed from the truth by a factor bounded only by a
function of $n$ (unfortunately, this function was of the order of magnitude
$n^{n/2}$).

In this paper we write down a formula (\eqref{iratio}) expressing the surface
area of $E$ in terms of an integral of a simple function over the sphere
$\mathbb{S}^{n-1}.$ This formula will be used to deduce a number of results:
\begin{enumerate}
\item The ratio of the surface area to the volume of $E$ (call
this ratio
  $\mathcal{R}(E)$) is a \emph{norm} on
  the vectors of inverse semi-axes. (Theorem \ref{isnorm}).
\item By a simple transformation (introduced for this purpose in \cite{riv1},
  though doubtlessly known for quite some time) $\mathcal{R}(E)$ can be
  expressed as a moment of a sum of independent Gaussian random
  variables; this transformation can be used to evaluate or estimate quite a
  number of related spherical integrals (see Section \ref{sphints}).
\item Quite  sharp bounds (\eqref{shbds}) on the ratio of
  $\mathcal{R}(E)$ to the $L^2$   norm of the vectors of inverses of
  semi-axes are derived.
\item We write down a very simple  asymptotic formula
  (Theorem \ref{asymp}) for the surface area of an ellipsoid of a very large
  dimension   with ``not too different'' axes. In particular, the formula
  holds if the ratio of the lengths of any two semiaxes is bounded by some
  fixed constant (Corollary \ref{asymp1}).
\item Finally, we give an identity relating the surface area of $E$ to a
  linear combination of Lauricella hypergeometric functions.
\end{enumerate}

\medskip\noindent
\textbf{Notation.} Let $(S, \mu)$ be a measure space with $\mu(S) < \infty.$
We will use the notation
\[
\fint_S f(x) \,d\mu \stackrel{\text{def}}{=} \dfrac{1}{\mu(S)} \int_S
f(x) d\mu.
\]
In addition, we shall denote the area of the unit sphere $\mathbb{S}^n$ by
$\omega_n$ and we shall denote the volume of the unit ball $\mathbb{B}^n$ by
$\kappa_n.$

\section{Cauchy's formula}
Let $K$ be a convex body in $\mathbb{E}^n.$ Let $u \in \mathbb{S}^{n-1}$ be
a unit vector, and let us define $V_u(K)$ to be the (unsigned)
$n-1$-dimensional volume of the orthogonal projection of $K$ in the direction
$U.$ \emph{Cauchy's formula} (see \cite[Chapter 13]{santalo}) then states that
%should be boxed
\begin{equation}
\label{cauchy}
\boxed{V_{n-1}(\partial K) = \dfrac{n-1}{\omega_{n-2}}
\int_{\mathbb{S}^{n-1}}V_u(K)\,d\sigma =
(n-1)\dfrac{\omega_{n-1}}{\omega_{n-2}} \fint_{\mathbb{S}^{n-1}}
  V_u(K)\,d\sigma,}
\end{equation}
where $d\sigma$ denotes the standard area element on the unit sphere.

In the case where $K=E$ is an ellipsoid, given by
\[
E = \{\mathbf{x}\in \mathbb{E}^n \quad | \quad \sum_{i=1}^n q_i^2 x_i^2 = \leq
1\}
\]
there are several ways of computing $V_u,$ one such way can be found in
\cite{connos}. In any event, the result is:
\begin{equation}
\label{conform}
V_u(E) = \kappa_{n-1} \dfrac{\sqrt{\left(\sum_{i=1}^n u_i^2
    q_i^2\right)}}{\prod_{i=1}^n q_i}.
\end{equation}
Since
\[
V_n(E) = \dfrac{\kappa_n}{\prod_{i=1}^n q_i},
\]
we can rewrite Cauchy's formula \eqref{cauchy} for $E$ in the form:
\begin{equation}
\label{iratio}
\boxed{\mathbb{R}(E) \stackrel{\text{def}}{=} \dfrac{V_{n-1}(\partial
    E)}{V_n(E)} =  n \fint_{\mathbb{S}^{n-1}} \sqrt{\sum_{i=1}^n u_i^2
    q_i^2}\,d\sigma,}
\end{equation}
where $\mathbb{R}(E)$ is the \emph{isoperimetric ratio} of $E.$
\begin{theorem}
\label{isnorm}
The ratio $\mathbb{R}(E)$ is a \emph{norm} on the vectors $q$ of lengths of
semiaxes ($q = (q_1, \dots, q_n).$)
\end{theorem}
\begin{proof}
The integrand in the formula \eqref{iratio} is a norm.
\end{proof}
\begin{corollary}
\label{pineq}
There exist constants $c_{n, p}, C_{n,p},$ such that
\[
c_{n, p} \| q \|_p \leq \mathbb{R}(q) \leq C_{n, p} \|q\|_p,
\]
where $\|q\|_p$ is the $L^p$ norm of $q.$
\end{corollary}
\begin{proof} Immediate (since all norms on a finite-dimensional Banach space
  are equivalent).
\end{proof}
In the sequel we will find reasonably good bounds on the constants $c_{n, p}$
and $C_{n, p},$ but for the moment observe that if $a_i = 1/q_i,\quad i=1,
\dots, n,$ then
\begin{equation}
\label{reverse}
\|q\|_p = \prod_{i=1}^n q_i \sigma_{n-1}^{1/p}(a_1^p, \dots,
a_n^p),
\end{equation}
where $\sigma_{n-1}$ is the $n-1$-st elementary symmetric function. In
particular, for $p=1,$ Corollary \ref{pineq} together with Eq. \eqref{reverse}
gives the estimate of \cite{riv2} (only for the $0$-th mean curvature integral
and with (for now) ineffective constants -- the latter part will be remedied
directly).
To exploit the formula \eqref{iratio} fully, we will need a digression on
computing spherical integrals.

\section{Spherical integrals}
\label{sphints}
In this section we will prove the following easy but very useful Theorem:
\begin{theorem}
\label{homothm}
Let $f(x_1, \dots, x_n)$ be a homogeneous function on $\mathbb{E}^n$ of degree
$d$ (in other words, $f(\lambda x_1, \dots, \lambda x_n) = \lambda^d f(x_1,
\dots, x_n)$.) Then
\[
\Gamma\left(\frac{n+d}2\right)
\fint_{\mathbb{S}^{n-1}} f d \sigma =
  \Gamma\left(\frac{n}{2}\right)\mathbb{E}\left(f(\mathbf{X}_1, \dots,
  \mathbf{X}_n)\right),\]
where $\mathbf{X}_1, \dots, \mathbf{X}_n$ are independent random variables
with probability density $e^{-x^2}.$
\end{theorem}
\begin{proof}
Let
\[E(f) = \fint_{\mathbb{S}^{n-1}} f(x)\,d\sigma,\]
and let $N(f)$ be defined as $\mathbb{E}(f(\mathbb{X}_1, \dots,
\mathbb{X_n})),$
where $\mathbb{X}_i$ is a Gaussian random variable with mean $0$ and variance
$1/2,$ (so with probability density $\mathfrak{n}(x) = e^{-x^2},$) and
$\mathbb{X}_1, \dots, \mathbb{X}_n$ are independent.
By definition,
\begin{equation}
\label{ndef}
N(f)(n) = c_n \int_{\mathbb{E}^n} \exp\left(-\sum_{i=1}^n x_i^2\right)
f(x_1, \dots, x_n) \,d x_1 \dots d x_n,
\end{equation}
where $c_n$ is such that
\begin{equation}
\label{norm}
c_n \int_{\mathbb{E}^n} \exp\left(-\sum_{i=1}^n x_i^2\right)\, d x_1 \dots d
x_n = 1.
\end{equation}
We can rewrite the expression \eqref{ndef} for $N(f)$ in polar coordinates as
follows (using the homogeneity of $f$):
\begin{equation}
\Nmin(n) = c_n \vol \mathbb{S}^{n-1}\int_0^\infty e^{-r^2} r^{n+d-1} E(f)
\,d r =  c_n E(f) \int_0^\infty e^{-r^2} r^{n+d-1}\, dr.
\end{equation}
 Since, by the substitution $u = r^2,$
\[\int_0^\infty e^{-r^2} r^{n+d-1}\, dr = \frac{1}{2}\int_0^\infty e^{-u}
u^{(n+d-2)/2} d u = \frac{1}{2}\Gamma\left(\frac{n+d}2\right).\]
and Eq. (\ref{norm}) can be rewritten in polar coordinates as
\[
1=c_n \vol \mathbb{S}^{n-1} \int_0^\infty r^{n-1} dr = \dfrac{c_n \vol
  \mathbb{S}^{n-1}}2 \Gamma\left(\frac{n}2\right),\]  we see that
\[\Gamma\left(\frac{n+d}2\right)E(f) = \Gamma\left(\frac{n}{2}\right)N(f).\]
\end{proof}
\begin{remark}
\label{gammarat}
In the sequel we will frequently be concerned with asymptotic results, so it
is useful to state the following asymptotic formula (which follows immediately
from Stirling's formula):
\begin{equation}
\label{gammaeq}
\lim_{x\rightarrow \infty} \dfrac{\Gamma(x+y)}{\Gamma(x)(x+y)^y} = 1.
\end{equation}
It follows that for large $n$ and fixed $d,$
\begin{equation}
\label{factorasymp}
\dfrac{\Gamma\left(\frac{n}2\right)}{\Gamma\left(\frac{n+d}2\right)} \sim
\left(\frac2{n+d}\right)^{d/2}.
\end{equation}
\end{remark}
\section{An explicit formula for the surface area}
The Theorem in the preceding section can be used to give explicit formulas
for the surface area of an ellipsoid (this formula will not be used in the
sequel, however). Specifically, in the book \cite{quad} there are formulas for
the moments of of random variables which are quadratic forms in Gaussian
random variables. We know that for our ellipsoid $E,$
\[\mathbb{R}(E) = n \fint_{\mathbb{S}^{n-1}} \sqrt{\sum_{i=1}^n u_i^2
  q_i^2}\,d\sigma = n
  \dfrac{\Gamma\left(\frac{n}2\right)}{\Gamma\left(\frac{n+1}2\right)}
\mathbb{E}\left(\sqrt{q_1^2 \mathbb{X}_1 +  \dotsb + q_n^2
  \mathbb{X}_n}\right),\] where $\mathbb{X}_i$ is a Gaussian with variance
  $1/2.$ The expectation in the last expression is the $1/2$-th moment of the
  quadratic form in Gaussian random variables, and so the results of
  \cite[p.~62]{quad} apply verbatim, so that we obtain:
\begin{equation}
\label{quadform1}
\mathbb{R}(E) = n
\dfrac{\Gamma\left(\frac{n}2\right)}{\Gamma\left(\frac{n+1}2\right)\Gamma\left(\frac12\right)}\sqrt{\alpha}
\int_0^\infty \dfrac{1}{\sqrt{z}} \sum_{j=1}^n \dfrac{q_j^2}{2(1+\alpha z
  q_j^2)} \left(\prod_{j=1}^n(1-q_j^2 z)\right)^{-1/2}\,dz;
\end{equation}
note that $\alpha$ in the above formula can be any positive number (as long as
 $|1-\alpha q_j^2| < 1,$ for all $j.$

This can also be expressed in terms of special functions.
 First, we need a definition:
\begin{definition}
Let $a, b_1, \dots, b_n, c, x_1, \dots, x_n$ be complex numbers, with
$|x_i| < 1,\quad i=1, \dots, n$, $\Re a > 0,$ $\Re(c-a) > 0.$ We then define
the \emph{Lauricella Hypergeometric Function} $F_D(a; b_1, \dots, b_n; c; x_1,
\dots, x_n)$ as follows:
\begin{multline}
F_D(a; b_1, \dots, b_n; c; x_1, \dots, x_n) = \\
\dfrac{\Gamma(c)}{\Gamma(a)\Gamma(c-a)} \int_0^1
u^{a-1}(1-u)^{c-a-1}\prod_{i=1}^n(1-ux_i)^{-b_i}\,d u.
\end{multline}
We also have the series expansion:
\begin{multline}F_D(a; b_1, \dots, b_n; c; x_1, \dots, x_n) = \\
\sum_{m_1 = 0}^\infty \cdots \sum_{m_n = 0}^\infty
\dfrac{(a)_{m_1 + \cdots + m_n} \prod_{i=1}^n (b_i)_{m_i}}{(c)_{m_1 + \cdots +
    m_n}} \prod_{i=1}^n \dfrac{x_i^{m_i}}{m_i!},
\end{multline}
valid whenever $|x_i| < 1, \forall i.$
\end{definition}
Now, we can write
\begin{multline}
\label{quadform2}
\mathbb{R}(E) = n
\dfrac{\Gamma^2\left(\frac{n}2\right)}{\Gamma^2\left(\frac{n+1}2\right)}
\sqrt{\alpha}\times\\
 \sum_{j=1}^n \frac{q_j^2}{2}
F_D\left(1/2; \eta_{1j}, \dotsc, \eta_{nj}; \frac{n+1}2; 1-\alpha q_1^2, \dotsc,
1-\alpha q_n^2\right),
\end{multline}
where $\eta_{ij} = 1/2 + \delta_{ij},$ and $\alpha$ is a positive parameter
satisfying $|1-\alpha q_j^2| < 1.$

\section{Laws of large numbers}
Many of the results in this section will require the following basic lemmas.
\begin{lemma}
\label{twomo}
Let $F_1, \dotsc, F_n, \dotsc$ be a sequence of probability distributions
whose first moments converge to $\mu$ and whose second moments converge to
$0.$ then $F_i$ converge to the Dirac delta function distribution centered on
$\mu.$
\end{lemma}
\begin{proof}
Follows immediately from Chebyshev's inequality.
\end{proof}
\begin{lemma}
\label{monotone}
Suppose the distributions $F_1, \dotsc, F_n, \dotsc$ converge to the
distribution $F,$ and the expectations of $|x|^\alpha$ with respect to $F_1,
\dots, F_n, \dotsc$ are bounded. Then the expectation of $|x|^\beta,$ $0
\leq \beta < \alpha$ converges to the expectation of $|x|^\beta$ with respect
to $F.$
\end{lemma}
\begin{proof}
See \cite[pp.~251-252]{feller2}.
\end{proof}
\begin{theorem}
\label{basiclaw}
Let $\mathbb{Y}_1, \dotsc, \mathbb{Y}_n, \dotsc$ be \emph{independent} random
variables with means
$0<\mu_1, \dotsc, \mu_n, \dotsc < \infty$ and variances $\sigma_1^2, \dotsc,
\sigma_n^2, \dotsc < \infty$ such that
\begin{equation}
\label{lovar}
\lim_{n \rightarrow \infty} \dfrac{\sum_{i=1}^n \sigma_i^2}{\left(\sum_{i=1}^n
  \mu_i \right)^2} = 0.
\end{equation}
Then
\[
\lim_{n\rightarrow \infty}\mathbb{E}\left(
\dfrac{\mathbb{Y}_1 + \dotsb + \mathbb{Y}_n}{\sum_{i=1}^n \mu_i}\right)^\alpha
= 1,\] for  $\alpha < 2.$
\end{theorem}
\begin{proof}
Consider the variable
\[
\mathbb{Z}_n = \dfrac{\sum_{i=1}^n \mathbb{Y}_i}{\sum_{i=1}^n \mu_i}.
\]
It is not hard to compute that
\[
\sigma^2(\mathbb{Z}_n) = \dfrac{\sum_{i=1}^n \sigma_i^2}{\left(\sum_{i=1}^n
  \mu_i \right)^2},
\]
while
\[\mu(\mathbb{Z}_n) = 1,\]
so by assumption \eqref{lovar} and Lemma \ref{twomo} $\mathbb{Z}_n$ converges
in distribution to the delta function centered at $1.$ The conclusion of the
Theorem then follows from Lemma \ref{monotone}.
\end{proof}
\begin{lemma}
\label{expcomp}
Let $\mathbb{X}$ be normal with mean $0$ and variance $1/2$ (so probability
density $e^{-x^2}/\sqrt{\pi}.$) Then
\[\mathbb{E}(|\mathbb{X}|^p) =
\dfrac{\Gamma\left(\frac{p+1}2\right)}{\sqrt{\pi}}.\]
\end{lemma}
\begin{proof}
\[\mathbb{E}(|\mathbb{X}|^p) = \frac{2}{\sqrt{\pi}}\int_0^\infty x^p
e^{-x^2}\,dx =
\frac{1}{\sqrt{\pi}}\int_0^\infty u^{(p-1)/2} e^{-u} =
\dfrac{\Gamma\left(\frac{p+1}2\right)}{\sqrt{\pi}}.\]
\end{proof}
\begin{theorem}
\label{normasymp}
\[
\fint_{\mathbb{S}^{n-1}}\|u\|_p \,d \sigma \sim
\dfrac{\Gamma\left(\frac{n}2\right)}{\Gamma\left(\frac{n+1}2\right)} \left(n
  \dfrac{\Gamma\left(\frac{p+1}2\right)}{\sqrt{\pi}}\right)^{\frac{1}{p}}.
\]
\end{theorem}
\begin{proof}
This follows immediately from the $1$-homogeneity of the $L^p$ norm, the
results of Section \ref{sphints}, Theorem \ref{basiclaw}, and Lemma
\ref{expcomp}.
\end{proof}
\subsection{Asymptotics of $\mathbb{R}(E).$}
\begin{theorem}
\label{asymp}
Let $q_1, \dotsc, q_n, \dotsc$ be a sequence of positive numbers such that
\[\lim_{n\rightarrow \infty} \dfrac{\sum_{i=1}^n q_i^4}{\left(\sum_{i=1}^n
  q_i^2 \right)^2} = 0.
\]
Let $E_n$ be the ellipsoid in $\mathbb{E}^n$ with major semiaxes $a_1 = 1/q_1,
\dotsc, a_n = 1/q_n.$ Then
\[
\lim_{n\rightarrow \infty}
  \dfrac{\Gamma\left(\frac{n+1}2\right)}{\Gamma\left(\frac{n}2\right)}
\dfrac{\mathbb{R}(E_n)}{n\sqrt{\frac12\sum_{i=1}^n q_i^2}} = 1.
\]
\end{theorem}
\begin{proof}
The Theorem follows immediately from Theorem \ref{basiclaw} and the results of
Section \ref{sphints}.
\end{proof}
\begin{corollary}
\label{asymp1}
Let $a_1, \dotsc, a_n, \dotsc$ be such that $0< c_1 \leq a_i/a_j \leq c_2 <
\infty,$ for any $i, j.$ Let $E_n$ be the ellipsoid with major semi-axes
$a_1, \dotsc, a_n.$ Then
\[
\lim_{n\rightarrow \infty}
  \dfrac{\Gamma\left(\frac{n+1}2\right)}{\Gamma\left(\frac{n}2\right)}
\dfrac{\mathbb{R}(E_n)}{n\sqrt{\frac12\sum_{i=1}^n \frac{1}{a_i^2}}} = 1.
\]
\end{corollary}
\begin{proof}
The quantities $q_1 = 1/a_1, \dotsc, q_n = 1/a_n, \dotsc$ clearly satisfy the
hypotheses of Theorem \ref{asymp}
\end{proof}
\section{General bounds on $\mathbb{R}(E)$}
We know that $\mathbb{R}(E)$ is a norm on the vector $\mathbf{q}=(q_1, \dotsc,
q_n)$  -- let us agree to write
\[
\|\mathbf{q}\|_\mathbb{R} \stackrel{\text{def}}{=}
\dfrac{\mathbb{R}(E)}{n}  = \fint_{\mathbb{S}^{n-1}}
\sqrt{\sum_{i=1}^n q_i^2 x_i^2} \,d \sigma.
\]
where $\mathbf{q}$ is the vector of inverses of the major semi-axes of $E.$
We know that for any $p > 0,$
\[ c_{n, p} \| q\|_p \leq \|q\|_\mathbb{R} \leq C_{n, p}\|q\|_p,\] for some
dimensional constants $c_{n, p}, C_{n, p}.$ In this section we
will give good ( estimates on the constants $c_{n,
  2}$ and $C_{n, 2}.$

To estimate $C_{n, 2}$ we will first show the following:
\begin{lemma}
\label{moments} Let $F(x)$ be a probability distribution, and let
\[M_a(F) \stackrel{\text{def}}{=} \mathbb{E}_F(|x|^a)\] denote the
absolute moments of $F$ (we will abuse notation in the sequel by
referring to the absolute moments of a random variable as well as
those of its distribution function). Further, let $0\leq \beta
\leq \alpha.$ Then
\[M_\beta(F) \leq 1 + M_\alpha(F).\]
\end{lemma}
\begin{proof}
\[M_\beta(F) = \int_{-\infty}^\infty |x|^\beta\,dF = \int_{-1}^1
|x|^\beta \,dF + \int_{-\infty}^{-1} |x|^\beta\,d F +
\int_1^\infty |x^\beta\,d F.
\]
We note that $\int_{-1}^1 |x|^\beta \,dF \leq \int_{-1}^1dF \leq
1,$ while $\int_{-\infty}^{-1} |x|^\beta\,d F \leq
\int_{-\infty}^0 |x|^\alpha \,d F,$ and similarly with the
integral from $1$ to $\infty.$ The assertion of the Lemma is then
immediate.
\end{proof}
\begin{theorem}\label{moments2}
Let $X_1, \dotsc, X_n$ be \emph{positive} i.~i.~d. random variables, with
finite mean $\mu.$ Let $a_1, \dotsc,
a_n$ be non-negative coefficients. Then
\[S=\mathbb{E}(\sqrt{\sum_{i=1}^n a_i X_i})\leq (1+\mu) \sqrt{\sum_{i=1}^n a_i}.\]
\end{theorem}
\begin{proof}
Let \[Y = \sum_{i=1}^n a_i X_i;\] it follows that $S =
M_{1/2}(Y).$ Let
\[Z = \dfrac{Y}{\sum_{i=1}^n a_i}.\] We see that
$M_1(Z) = \mu.$ By Lemma \ref{moments},
it follows that
\[M_{1/2}(Z) \leq \mu + 1,\]
and thus \[S = M_{1/2}(Y) \leq
(1+\mu) \sqrt{\sum_{i=1}^n a_i}.\]
\end{proof}
\begin{corollary}
\label{qbound}
\[C_{n, 2} \leq
\dfrac{3\Gamma\left(\frac{n}2\right)}{2\Gamma\left(\frac{n+1}2\right)}
.\]
\end{corollary}
\begin{proof}
Replace $a_i$ by $q_i^2$ and make the variables $X_i$ squares of centered
Gaussians  with
variance $1/2$ in the statement of Theorem \ref{moments2}, and
apply the results of Section \ref{sphints} to get the quotient of
$\Gamma$ values.
\end{proof}
To estimate $c_{n, 2}$ we will first note:
\begin{lemma}
\label{convsqrt}
Let $\lambda_1, \dotsc, \lambda_n \geq 0,$ and let $\lambda_1 + \dotsb
+ \lambda_n = 1.$ 
Then
\[
\sqrt{\sum_{i=1}^n \lambda_i x_i^2} \geq \sum_{i=1}^n \lambda_i |x_i|.
\]
\end{lemma}
\begin{proof}
Concavity of square root.
\end{proof}
\begin{corollary}
\label{preqbound}
\[
\fint_{\mathbb{S}^{n-1}} \sqrt{\sum_{i=1}^n q_i^2 u_i^2} \,d \sigma
  \geq \fint_{\mathbb{S}^{n-1}} |x_1| \,d\sigma 
\sqrt{\sum_{i=1}^n  q_i^2}
.\]
\end{corollary}
\begin{proof}
Write
\[ \sum_{i=1}^n q_i^2 u_i^2 = \sum_{i=1}^n q_i^2
\sum_{i=1}^n \frac{q_i^2}{\sum_{j=1}^n q_j^2} x_i^2,\] then use the
symmetry to note that \[\fint_{\mathbb{S}^{n-1}} x_i \,d\sigma =
=\fint_{\mathbb{S}^{n-1}} x_j \,d\sigma\] for any $i, j.$
\end{proof}
\begin{corollary}
\label{qqbound}
\[c_{n, 2} \geq
\dfrac{\Gamma\left(\frac{n}2\right)}{\sqrt{\pi}\Gamma\left(\frac{n+1}2\right)}.\] 
\end{corollary}
\begin{proof}
Immediate from Corollary \ref{preqbound}
\end{proof}
Combining Corollary \ref{qbound} and Corollary  \ref{qqbound} we get
\begin{equation}
\label{shbds} 
\dfrac{3\Gamma\left(\frac{n}2\right)}{2\Gamma\left(\frac{n+1}2\right)}
\geq
\dfrac{\|q\|_{\mathbb{R}}}{\|q\|} \geq
\dfrac{\Gamma\left(\frac{n}2\right)}{\sqrt{\pi}\Gamma\left(\frac{n+1}2\right)}.
\end{equation}
Note that the ratio between the right and
the left hand sides of the inequality \eqref{shbds} stays bounded,
 so
this inequality is sharp to within a constant factor. 
It is fairly clear that the constants are not sharp, however.

\bibliographystyle{amsplain}

\end{document}